\documentclass[12pt]{article}
\usepackage{amsfonts, amssymb, amsmath}
\usepackage{amsthm}
\usepackage[all]{xy}
\usepackage{url}
\usepackage{graphicx}
\usepackage{listings}
\newcommand{\zz}{\mathbb{Z}}

\newcommand{\e}{\equiv}
\newcommand{\el}{\lambda}
\newcommand{\x}{\cdot}

\usepackage[parfill]{parskip}

\begin{document}
\begin{center}
\textbf{\Large A note on some properties of the $\lambda$-Polynomial} \\
\vspace{5mm}
\begin{tabular}{c} David Bodiu\\
\multicolumn{1}{c}{Trinity College Dublin} \\
\verb+bodiud@tcd.ie+ 
\end{tabular} \end{center}

\begin{abstract}
\noindent The expression $a^n + b^n$ can be factored as $(a+b)(a^{n-1} - a^{n-2} b + a^{n-3} b^2 - ... + b^{n-1})$ when $n$ is an odd integer greater than one. This paper focuses on proving a few properties of the longer factor above, which we call $\el_n(a,b)$. One such property is that the primes which divide $\el_n(a,b)$ satsify $p \ge n$, if $a,b$ are coprime integers and $n$ is an odd prime. \end{abstract}

\section{Preliminary lemmas}

\textbf{Lemma 1.1} Let $a,b \in \zz$ and let $d \in \zz$ divide $a$ but not $b$. Then $d$ does not divide $a+b$.\\

\textit{Proof.} Suppose otherwise. We can then write $a = kd$ and $a+b = hd$ for some $k,h \in \zz$. As a consequence $b = hd - kd = d(h-k)$ which implies that $b$ is divisible by $d$, a contradiction.
\null \hfill{$\blacksquare$} \\

\textbf{Lemma 1.2} Let $a, b \in \zz$ be coprime. Then $ab$ and $a+b$ are coprime.  \\

\textit{Proof.} Notice that $ab = \pm 1$ only if $a,b \in \{1,-1\}$, in which case $ab$ and $a+b$ are coprime. Otherwise suppose $p$ is a prime which divides $ab$. Then $p$ divides either $a$ or $b$, but not both. It therefore follows by Lemma 1.1 that $p$ cannot divide $a+b$.
\null \hfill{$\blacksquare$} \\

\textbf{Corollary 1.3} Let $a, b \in \zz$ be coprime. Then $(ab)^{m_1}$ and $(a+b)^{m_2}$ are coprime for any positive integers $m_1$ and $m_2$. \\

\textit{Proof.} By Lemma 1.2 it follows that $ab$ and $a+b$ are coprime. But raising any of these to some positive power does not change their prime factors, which means that they remain coprime.
\null \hfill{$\blacksquare$} \\

\section{The $\el$-Polynomial}

Suppose that $n$ is an odd integer greater than one and $a,b \in \zz$. In this section we wish to study the longer factor of $a^n+b^n=(a+b)(a^{n-1} - a^{n-2} b + a^{n-3} b^2 - ... + b^{n-1})$ which we will call $\el_n(a,b)$. Notice that when $b \neq -a$ we can rearrange the above equation to obtain

$$\el_n(a,b) = \frac{a^n + b^n}{a+b}.$$

In the case that $b = -a$ the original expression gives $\el_n(a,-a) = n a^{n-1}$. This provides us with two equivalent definitions for $\el_n(a,b)$: \\

\textbf{Definition 2.1} Let $a,b \in \zz$ and $n$ be an odd integer greater than one. We define 

$$\el_n(a,b) = \sum_{i=0}^{n-1} (-1)^i a^{n-1-i} b^i$$

or

$$\el_n(a,b) =
\begin{cases}
\cfrac{a^n + b^n}{a+b} &\text{if } b \neq -a; \\
na^{n-1} &\text{if } b = -a.
\end{cases}$$

The second definition is much more compact and just nicer work with in general. \\

\textbf{Lemma 2.2} (Non-negativity) Let $a,b \in \zz$ and $n$ be an odd integer greater than one. Then $\el_n(a,b)$ is a non-negative integer. \\
 
 \textit{Proof.} In the case $b \neq -a$ we have $\el_n(a,b) = \frac{a^n +b^n}{a+b}$. Because $n$ is odd the numerator and denominator will always have the same sign, ensuring positivity. In the case $b = -a$ we have $\el_n(a,b) = n a^{n-1}$. Since $n$ is positive and $a^{n-1}$ is non-negative we must have that their product is also non-negative.
 \null \hfill{$\blacksquare$} \\
 
 \textbf{Remark 2.3} Notice that $\el_n(a,b)$ is actually a positive integer if we add the condition that at least one of $a$ and $b$ are nonzero. Note also that this condition is implied if $a$ and $b$ are coprime. \\

\textbf{Lemma 2.4} (Symmetry) Let $a,b \in \zz$ and let $n$ be an odd integer greater than one. Then $\el_n(a,b) = \el_n(b,a)$. \\

In the case $b \neq -a$ this is quite plain. In the case $b = -a$ we have that $\el_n(a,-a) = na^{n-1}$ is the same as $\el_n(-a,a) = n(-a)^{n-1}$ because $n-1$ is even.
 \null \hfill{$\blacksquare$} \\

  \textbf{Lemma 2.5} (Homogeneity) Let $a,b,d \in \zz$ and let $n$ be an odd integer greater than one. Then $\el_n(da,db) = d^{n-1} \cdot \el_n(a,b)$.\\
  
\textit{Proof.} This is quite easy to see both when $b \neq -a$ and $b = -a$.
\null \hfill{$\blacksquare$} \\
  
\textbf{Corollary 2.6} (Sign Symmetry) Let $a,b \in \zz$ and let $n$ be an odd integer greater than one. Then $\el_n(-a,b) = \el_n(a,-b)$. \\

\textit{Proof.} By Lemma 2.5

$$\el_n(-a,b) = \el_n(-1 \cdot a, -1 \cdot (-b)) = (-1)^{n-1} \el_n(a,-b) = \el_n(a,-b).$$
\null \hfill{$\blacksquare$} \\

\textbf{Corollary 2.7} (Eveness) Let $a,b \in \zz$ and let $n$ be an odd integer greater than one. Then $\el_n(-a,-b) = \el_n(a,b)$. \\

\textit{Proof.} By Lemma 2.5

$$\el_n(-a,-b) =\el_n(-1 \cdot a, -1 \cdot b) = (-1)^{n-1} \el_n(a,b) = \el_n(a,b).$$
  \null \hfill{$\blacksquare$} \\
  
The next few results are more selective, requiring that $a$ and $b$ be coprime integers or that $n$ be an odd prime. \\

\textbf{Lemma 2.8} Let $a, b \in \zz$ be coprime and $n$ be an odd integer greater than one. Then $\el_n(a,b)$ is coprime to $a$ and $b$. \\

\textit{Proof.} By using the suitable values for Lemma 1.2, we see that $a^n + b^n$ is coprime to both $a$ and $b$. Since $a^n + b^n = (a+b) \x \el_n(a,b)$. It follows that $\el_n(a,b)$ must also be coprime to $a$ and $b$.
  \null \hfill{$\blacksquare$} \\

\textbf{Proposition 2.9} Let $a, b \in \zz$ be coprime and $n$ be an odd integer greater than one. Then the common divisors of $a+b$ and $ \el_n(a,b)$ must divide $n$. \\

\textit{Proof.} We consider the expression

$$\el_n(a,b) \pm n a^{\frac{n-1}{2}}b^{\frac{n-1}{2}}$$

with ``$-$'' if $\frac{n-1}{2}$ is even and ``$+$" if $\frac{n-1}{2}$ is odd. If we can prove that the above expression is divisible by $a+b$, then we'll be able to rearrange some things and prove the Proposition. We therefore view the expression modulo $a+b$. Since $b \e -a \pmod{a+b}$ we get

$$\el_n(a,b) \pm n a^{\frac{n-1}{2}}b^{\frac{n-1}{2}} \e na^{n-1} \pm na^{\frac{n-1}{2}}(-a)^{\frac{n-1}{2}} \pmod{a+b}.$$

If $\frac{n-1}{2}$ is even we have that

$$\el_n(a,b) - n a^{\frac{n-1}{2}}b^{\frac{n-1}{2}} \e n a^{n-1} - n a^{n-1} \e 0 \pmod{a+b}$$

otherwise if $\frac{n-1}{2}$ is odd

$$\el_n(a,b) + n a^{\frac{n-1}{2}}b^{\frac{n-1}{2}} \e n a^{n-1} - n a^{n-1} \e 0 \pmod{a+b}.$$

Consequently, $\el_n(a,b) + k(a+b) = \mp na^\frac{n-1}{2}b^\frac{n-1}{2}$ for some $k \in \zz$. By Lemma 2.8, $\el_n(a,b)$ is coprime to $a$ and $b$, therefore $k \neq 0$. Suppose then that $d$ divides both $a+b$ and $\el_n(a,b)$. It follows by the equation above that $d$ also divides $\mp na^\frac{n-1}{2}b^\frac{n-1}{2}$. However, applying Corollary 1.3 with $m_1 = \frac{n-1}{2}$ and $m_2 = 1$ it follows that $d$ must divide $n$.
 \null \hfill{$\blacksquare$} \\

\textbf{Proposition 2.10} Let $a,b \in \zz$ and $n$ be an odd prime.  Then $a+b$ is divisible by $n$ if and only if $\el_n(a,b)$ is divisible by $n$. \\

\textit{Proof.} Suppose that $a+b$ is divisible by $n$. Then $b \e -a$ (mod $n$), and as a consequence $\el_n(a,b) \e n a^{n-1} \e 0$ (mod $n$). Suppose on the other hand that $\el_n(a,b)$ is divisible by $n$. Since $(a+b) \x \el_n(a,b) = a^n +b^n$ we obtain that $a^n + b^n \e 0$ (mod $n$). Finally, applying Fermat's Little Theorem gives $a+b \e 0$ (mod $n$).
 \null \hfill{$\blacksquare$} \\

 The following result is taken from \cite{gauss1986disquisitiones} and is provided without proof: \\

\textbf{Lemma 2.11} Let $p$ be a prime and $m$ be a positive integer. If $m$ is coprime to  $p-1$, then every element of $\zz / p \zz$ has a unique $m$-th root.
 \null \hfill{$\blacksquare$} \\

\textbf{Theorem 2.12} Let $a,b \in \zz$ be coprime and let $n$ be an odd prime. Suppose that $p \neq n$ is a prime which divides $\el_n(a,b)$. Then $p \e 1$ (mod $n$). \\

\textit{Proof.} We once again use the fact that $\el_n(a,b) \cdot (a+b) = a^n +b^n$, from which we obtain $a^n + b^n \e 0$ (mod $p$), or rearranged: $a^n \e (-b)^n$ (mod $p$). Now by Lemma 2.11, if $n$ does not divide $p-1$, we are able to take $n$-th roots modulo $p$. Suppose that this is the case. Then, $a \e -b$ (mod $p$), that is, $a+b \e 0$ (mod $p$). However, by Proposition 2.9, this is a contradiction. We must therefore have that $n$ divides $p-1$, in other words, $p-1 \e 0$ (mod $n$), which completes the proof.
\null \hfill{$\blacksquare$} \\

\textbf{Remark 2.13} Notice that if $a$ and $b$ are not coprime, by Lemma 2.5 the common prime factors of $a$ and $b$ will form ``clusters" of the form $p^{n-1}$ in $\el_n(a,b)$. These ``clusters" are congruent to one modulo $n$ by Fermat's Little Theorem, as long as $p \neq n$. The other prime factors of $\el_n(a,b)$ will be each individually congruent to one modulo $n$ by Theorem 2.12, as long as they are not equal to $n$.\\

\textbf{Corollary 2.14} Let $a,b \in \zz$ be coprime and let $n$ be an odd prime. If $p \neq n$ is a prime which divides $\el_n(a,b)$, then $p > n$. \\

\textit{Proof.} Suppose not. Then there exists some prime factor $p$ of $\el_n(a,b)$ which satisfies $p<n$. By Theorem 2.12 we also have that $p \e 1$ (mod $n$). This implies that $p = 1$, a contradiction.
\null \hfill{$\blacksquare$} \\

\textbf{Corollary 2.15} Let $a,b \in \zz$ be coprime and let $n$ be an odd prime. Then $\el_n(a,b)$ is an odd integer. \\

\textit{Proof.} Corollary 2.14.
\null \hfill{$\blacksquare$} \\

\textbf{Corollary 2.16} Let $a,b \in \zz$ be coprime and let $n$ be an odd prime. If $n$ does not divide $\el_n(a,b)$, then $\el_n(a,b) \e 1$ (mod $n$). \\

\textit{Proof.} The result follows from combining Lemma 2.2 with Theorem 2.12.
\null \hfill{$\blacksquare$} \\

\textbf{Lemma 2.17} (Generalization of Corollary 2.16) Let $a,b \in \zz$, not necessarily coprime, and let $n$ be an odd prime. If $n$ does not divide $\el_n(a,b)$, then $\el_n(a,b) \e 1$ (mod $n$). \\

\textit{Proof.} This can be deduced directly from Lemma 2.2 and Remark 2.13. Alternatively, we have that $\el_n(a,b) \cdot (a+b) = a^n + b^n$. Using Fermat's Little Theorem we obtain $\el_n(a,b) \cdot (a+b) \e a+b$ (mod $n$). By Proposition 2.10, since $n$ does not divide $\el_n(a,b)$ it does not divide $a+b$ either. As such, we can multiply both sides of the congruence above by $(a+b)^{-1}$ to get $\el_n(a,b) \e 1$ (mod $n$).
\null \hfill{$\blacksquare$} \\

\section*{Acknowledgements}
The author is grateful to Nicolas Mascot for his helpful comments on this paper.

\bibliographystyle{plain}
\bibliography{References}

\end{document}